\documentclass[12pt]{amsart}
\usepackage{amsmath,amsfonts,amssymb,amsthm}

\theoremstyle{definition}

\theoremstyle{remark}

\def\m1{^{-1}}
\def\ov1{\overline}

\title[Wreath products in modular group algebras $\dots$]
{Wreath products in modular group algebras of some finite 2-groups}
\author{Alexander Konovalov}
\dedicatory{Dedicated to Professor A.A.~Bovdi on his 70th birthday}
\address{
School of Computer Science, University of St Andrews, \newline
\indent North Haugh, St Andrews, Fife, KY16 9SX, Scotland}

\email{konovalov@member.ams.org}


\date{}

\begin{document}

\begin{abstract}
Let $K$ be field of characteristic 2 and let $G$ be a finite
non-abelian 2-group with the cyclic derived subgroup $G'$, and
there exists a central element $z$ of order 2 in $Z(G) \backslash
G'$. We prove that the unit group of the group algebra $KG$
possesses a section isomorphic to the wreath product of a group of
order 2 with the derived subgroup of the group $G$, giving for
such groups a positive answer to the question of A.~Shalev.
\end{abstract}

\maketitle

\section{Introduction}
\label{intro}

Let $p$ be a prime number, $G$ be a finite $p$-group and $K$ be a
field of characteristic \nolinebreak $p$. Denote by $I(KG)$ the 
augmentation ideal of the modular group algebra $KG$. The group 
of normalized units $V(KG)$ consists of all elements of the form 
$1+x$, where $x$ belongs to $I(KG)$.

An interest to the structure of the normalized unit group raised 
a number of questions about kinds of wreath products that may be 
involved into $V(KG)$ as a subgroup of as a section, i.e. as a 
factor-group of a certain subgroup of $V(KG)$.

The first result was obtained in \cite{Col-Pas} by D.~Coleman and 
D.~Passman, who proved that for a non-abelian finite $p$-group $G$ 
a wreath product of two groups of order $p$ is involved into $V(KG)$. 
Later it was generalized by A.~Bovdi in \cite{Bovdi3}. Among other 
related results it is worth to mention \cite{Mann,Mann-Sh,Sh6}.
C.~Bagi\'{n}ski in \cite{Ba} described all $p$-groups, for which 
$V(KG)$ does not contain a subgroup isomorphic to the wreath  
product of two groups of order $p$ for the case of odd $p$. 
Using results and methods of \cite{Bovdi-Kovacs}, the case of $p=2$ 
was investigated in \cite{Bovdi-Dokuchaev} by V.~Bovdi and M.~Dokuchaev.

In \cite{Shalev1} A.~Shalev, motivated by problem of determining
of the nilpotency class of $V(KG)$, formulated the question whether 
$V(KG)$ possesses a section isomorphic to the wreath product of a 
cyclic group of order $p$ and the derived subgroup of $G$. In
\cite{Shalev2} he proved that this is true for the case of an odd
$p$ and a cyclic derived subgroup of $G$.

Besides the importance of wreath products with large nilpotency
class in the investigation of the unit group, it appears that there are
certain connections between wreath products and Lie nilpotency indices 
of group algebras (see \cite{BJS,Bovdi-Spinelli}).

When $p=2$ and $G$ is a 2-group of almost maximal class, in
\cite{Konovalov3} and \cite{Konovalov4} the author
constructed a section of $V(KG)$ isomorphic to the wreath product
of a group of order 2 and the derived subgroup of $G$.

The aim of the present short note is to publish one observation,
that allows to confirm the conjecture of A.~Shalev for another 
class of 2-groups. 

\vskip 0.2cm
{\bf Theorem.} Let $K$ be field of characteristic 2 and let $G$ be
a finite non-abelian 2-group with the cyclic derived subgroup
$G'$, and there exists a central element $z$ of order 2 in
$Z(G) \backslash G'$. Then the wreath product of the cyclic group
of order 2 and the derived subgroup of $G$ is involved into the
$V(KG)$.
\vskip 0.2cm

\section{Proof of the theorem}
\label{wreath}

{\bf Proof.} Let $K$ be field of characteristic 2 and let $G$ be a
finite $p$-group of order $p^{n}$ with the cyclic derived
subgroup $G'= \langle (b,a) \rangle$ of order $2^{s}$, and there
exists an element $z$ of order 2 in $Z(G) \setminus G'$.
Let the order of $a$ is $p^k$ and $(b,a^{p^{s}})=1$.
Consider an element $h=1+b(1+z)$. Then it is easy to check that
\[
\begin{split}
h^{a} & =1+b^{a}(1+z)=1+b(b,a)(1+z), \\
h^{a^{2}} & =1+b^{a^{2}}(1+z)=1+b(b,a^{2})(1+z), \\
\cdots \\
h^{a^{p^{s}}-1} & =1+b^{a^{p^{s}}-1}(1+z)=1+b(b,a^{p^{s}-1})(1+z), \\
h^{a^{p^{s}}} & =1+b^{a^{p^{s}}}(1+z)=h. \\
\end{split}
\]
Note that $h$ is of order 2 as well as all of its conjugates, and
they generate an elementary abelian group $ \langle h, h^{a},
h^{a^{2}}, \dots, h^{a^{p^{s}}-1} \rangle $. It is easy to see
that this group is actually the direct product 
$$ X = \langle h
\rangle \times \langle h^{a} \rangle \times \langle h^{a^{2}}
\rangle \times \cdots \times \langle h^{a^{p^{s}}-1} \rangle, $$
since if we multiply its elements, we will obtain an element of
the form
$$ 1+b((b,a^{i_1})+(b,a^{i_2})+ \cdots + (b,a^{i_t}))(1+z)$$
which is not equal to 1. Now the wreath product could be obtained
if we factorize the semidirect product of $X$ and $\langle a
\rangle$ over $\langle a^{p^{s}} \rangle$, so the theorem is proved.

{\bf Remark.} Note that the family of groups from the conditions
of the theorem extends the result of \cite{Konovalov4}
significantly. For example, using the Small Groups Library of the
GAP system \cite{GAP}, we can find that the number of such groups
of order $2^n$ for $n=4,5,6,7,8,9$ is is accordingly 4, 20, 72,
231, 662 and 1750.

\nocite{*}
\bibliographystyle{abbrv}
\bibliography{wreath}

\end{document}